\documentclass{elsart}
\usepackage{graphics}
\usepackage{graphicx}
\usepackage{epsfig}
\usepackage{amssymb}
\usepackage{amsmath}
\usepackage{multirow}
\usepackage{theorem,amsfonts,psfig}
\usepackage{verbatim}
\usepackage{subfigure}
\usepackage{algorithm}
\usepackage{algpseudocode}
\usepackage{pifont}
\usepackage{lscape}
\usepackage{tikz,tkz-tab}
\tikzstyle{diam}=[rectangle,draw,fill=green!30]

\newtheorem{definition}{Definition}[section]
\begin{document}
\begin{frontmatter}
\title{Spherical harmonics entropy for optimal 3D modeling}
\author{Malika Jallouli}
\ead{jallouli.malika@yahoo.fr}
\address{Department of Electric Engineering, National Engineering School of Sousse, Tunisia.}
\author{Wafa Belhadj Khalifa}
\ead{wafa.bhk@gmail.com}
\address{Department of Electric Engineering, National Engineering School of Sousse, Tunisia.}
\author{Anouar Ben Mabrouk\thanksref{label2}}
\ead{anouar.benmabrouk@fsm.rnu.tn}
\address{Research Unit of Algebra, Number Theory and Nonlinear Analysis UR11ES50, Department of Mathematics, Faculty of Sciences, 5019 Monastir. Tunisia}
\thanks[label2]{Department of Mathematics, HIgher Institute of Applied Mathematics and Informatics, Street of Assad Ibn Al-Fourat, Kairouan 3100, Tunisia.}
\author{Mohamed Ali Mahjoub}
\ead{medalimahjoub@gmail.com}
\address{Department of Electric Engineering, National Engineering School of Sousse, Tunisia.}
\maketitle
\begin{abstract}
3D image processing constitutes nowadays a challenging topic in many scientific fields such as medicine, computational physics and informatics. Therefore, development of suitable tools that guaranty a best treatment is a necessity. Spherical shapes are a big class of 3D images whom processing necessitates adoptable tools. This encourages researchers to develop shperical wavelets and spherical harmonics as special mathematical bases able for 3D spherical shapes. The present work lies in the whole topic of 3D image processing with the special spherical harmonics bases. A spherical harmonics based approach is proposed for the reconstruction of images provided with spherical harmonics Shannon-type entropy to evaluate the order/disorder of the reconstructed image. Efficiency and accuracy of the approach is demonstrated by a simulation study on several spherical models.
\end{abstract}
\end{frontmatter}
\section{Introduction}
Many signals/images are naturally contained or supported and thus parametrized over a sphere. Concrete examples may be obtained from astronomy, computer graphics, medical imaging, ... etc. Therefore, an efficient representation of spherical signals/images is of great importance. The most important interest may be particularly related to the ability of mathematical/physical bases to approximate such objects accurately with minimal storage costs, and the possibility of obtaining computationally efficient algorithms to process them.

One motivation behind the present work is the introduction of such an efficient representation of spherical signals already with a suitable 'best' measure of the efficiency, such as entropy.

In recent decades there has been a great growth in 3D image theory as wll as applications. This leads 3D image processing to constitute nowadays a challenging topic in both theoretical and applied fields. Mathematicians, theoretical informaticians and theoretical physicists are developping suitable tools for best representing these 3D shapes which in turn will be applied in medicine, computational physics, informatics, ... etc. Spherical shapes are a big class of 3D images whom processing necessitates adoptable tools. This encourages researchers to develop shperical wavelets and spherical harmonics as the most efficient special mathematical bases able for 3D spherical shapes nowadays. The present work aims in a first phase to review the concepts of these topics to be applied later. Next, as a natural phase in the application of any theoretical tool is the checking out of its efficiency. We thus propose in the present work to 'introduce' some special form of entropy to attain this object.

Recall that, generally, the concept of entropy in information theory, image processing, dynamical systems, mathematics, dimension theory is not receent. It had origins in thermodynamics. It is considered as a measure of the system’s disorder that is directly effected by the reversibility/irreversibility of phase changement such as heat.

Entropy continues to fascinate researchers, engineers, mathematicians, chemists, computer scientists, physicists. In \cite{philBroadbridge}, the scope of entropy as a diagnostic tool in higher-order partial differential equations has been investigated. In \cite{Muller} introductory concepts and general physical problems on entropy have been conducted with a focus on the eventual link energy/entropy. In \cite{Muller2}, a theory and proofs are presented for a version of irreversible modeling. In \cite{Robinson} a mathematical framework has been developed provided with a study of the link between the entropy in information theory and entropy in dynamical systems and in statistical physics. In \cite{Ruggeri}, further developments on entropy have been investigated in the framework of hyperbolic conservation laws in continuum mechanics. Alternative formulations of irreversible thermodynamics are given in \cite{Attard} and \cite{Fradkov}. In \cite{Attard}, an entropy variant has been introduced as maximum quantity allowing to estimate an optimal state of a non-equilibrium system. In \cite{Fradkov} the speed gradient model of transient dynamics has been derived for finitely many systems. 

More recently, the concept of entropy has been generalized, extended and more and more investigated in many domains such as computer information theory, topological entropy, metric entropy of Kolmogorov-Sinai and dynamical systems in mathematics. 

In signal analysis (1D image), entropy is somehow well known and applied quite widely. In \cite{Nicholson2}, the degree of disorder in earthquake signals is measured using the concept of entropy. A correlation between earthquake entropy and tectonic regime has been discovered. A temporal definition of entropy has been applied in \cite{TelescaL} and \cite{TelescaL2} to study seismicity in different parts of the earth.

In 2D image, the concept of entropy is somehow new. It is introduced recently in the analysis of 2D images such as segmentation in \cite{BWanga} using maximum entropy thresholding. In \cite{OriettaN}, algorithms have been proposed to analyze the irregularity on spatial scales, which extends the 1D multi-scale entropy to the 2D case. In \cite{PereP}, the measure of entropy has been applied to derive and automatic localization of positions. 

However, in the most cited litterature, the estimation of entropy is based on the probabilistic/statistical classical formulations which is based on the evaluation of the probability of presence of and/or the probability anf frequency of ordering a studied system. Recall that in mathematical physics theory, the entropy is a mathematical function based on probabilistic statistical theory that intuitively corresponds to the amount of information contained or delivered by an information source.

The concept of wavelet based entropy still remains recent and the studies investigating such a concept are few. This motivates the development of more works about it. In \cite{OriettaN} a wavelet-based approach to measure Shannon's entropy in the context of spatial point models has been conducted. The spatial heterogeneity and complexity of spatial point models are analyzed using the entropy of multiscale anisotropic wavelets.

The present paper is subscribed in the whole topic of investigating entropy concept by means of wavelet theory. We aim towfolds. Firstly, we aim to extend such a concept to the 3D case and secondly to apply a special image processing functional analysis issued from spherical wavelets consisting of spherical harmonics. A correlation-type relation between entropy and sperical harmonics modeling is proposed. In fact, we intend to define, in a precise and automatic way, the final reconstruction order of spherical harmonics models that best represents initial images to be reconstructed. In this case, a new entropy based on the spherical harmonics decomposition is proposed and applied on 3D surfaces. 

The present paper will be thuerefore organized as follows. Spherical harmonics modeling is reviwed in section 2. Then, the concept of spherical hamonics entropy is introduced in section 3 with a breif reminder of Shannon's classical entropy. Experimentations are next conducted in order to prove the efficiency of the proposed theory. We conclude afterward.
\section{Harmonic analysis on the sphere: Spherical Harmonics}
The key point of harmonic analysis on the sphere is that it provides a type of analysis that lives entirely on the sphere, the signal/image and the analyzing basis. Mathematically speaking, harmonic analysis on the sphere is not new, and its development is widely investigated since the begining of the 1900'th century. Harmonic analysis on the sphere is the natural extension of Fourier series,
which studies the expressibility of functions and generalized functions as sums
of the fundamental exponential functions. The exponential functions are
simpler functions, and are both eigenfunctions of the translation-invariant
differential operator, and group homomorphisms too. Here also, spherical
harmonics are simple and eigenfunctions of some differential operators. More backgrounds may be found in \cite{Antoine-Carrette-Murenzi-Piette}, \cite{Antoine-Demanet-Jacques-Vandergheynst}, \cite{Antoine-Murenzi-Vandergheynst}, \cite{Arfaouietal2}, \cite{Bulow-Daniilidis}, \cite{Guan}.
\subsection{Classical harmonic analysis on the sphere}
In 3D image processing, harmonic analysis starts from the sphere $S^2\subset\mathbb{R}^3$ and analyzes the $L^2(S^2)$-functions by means of well adapted functional bases. The most known beses are the so-called spherical harmonics which are the everywhere regular eigenfunctions of the Laplacian operator on the sphere $S^2$, 
$$
\medskip\Delta=\displaystyle\frac{1}{r^2}\displaystyle\frac{\partial}{\partial r}\left(r^2\displaystyle\frac{\partial}{\partial r}\right)+\frac{1}{r^2\sin\theta}\displaystyle\frac{\partial}{\partial\theta}\left(\sin\theta\displaystyle\frac{\partial}{\partial\theta}\right)+\frac{1}{r^2\sin^2\theta}\displaystyle\frac{\partial^2}{\partial\varphi^2}=0.
$$
Spherical harmonics are usually charterized by two parameters $l\in\mathbb{N}$ and $m\in\mathbb{Z}$ such that $|m|\leq l$. Denote $Y_{l,m}$ the associated spherical harmonic to $(l,m)$. It is expressed as 

\begin{equation}\label{SHbasis}
Y_{l,m}=K_{l,m}P_{l,m}(cos\theta)e^{im\varphi}
\end{equation}

where $K_{l,m}$ is the normalization constant
$$
K_{l,m}=(-1)^m\sqrt{\frac{(2l+1)(l-m)!}{4\pi(l+m)!}}.
$$
$P_{l,m}$ is the well known Legendre polynomial of degree $l$ and order $m$,
$$
P_{l,m}(x)=\frac{(-1)^m}{2^ll!}(1-x^2)^{\frac{m}{2}}\frac{d^{l+m}}{dx^{l+m}}(x^2-1)^l
$$
The set $(Y_{l,m})$ constitutes an orthonormal basis of the space $L^2(S^2)$ of all images with finite variance or energy. The basis of spherical harmonics ($Y_{l,m}$) has several mathematical properties useful for modeling and pattern recognition such as
\begin{description}
	\item[i.] $S^2$ being a compact group, the Fourier transform is therefore represented by the Fourier coefficients with respect to the associated Legendre base.
	\item[ii.] Spherical harmonics form a complete set on the surface of the unit sphere.
\end{description}
These two properties make it possible to deduce the reconstruction formula of the surface to be modeled. Mathematically speaking, any function $f\in L^2(S^2)$ can be decomposed in a series called the spherical harmonics' series decomposition,
\begin{equation}\label{decompositionshseries}
f(\theta,\varphi)=\sum_{l=0}^{\infty}\sum_{m=-l}^l S_{l,m}Y_{l,m}(\theta,\varphi)
\end{equation}
where for $(l,m)$ fixed, $S_{l,m}$ is the harmonic coefficient in the degree $l$ and the order $m$ evaluated by analogy to the 2D Fourier coefficients,
\begin{equation}\label{shcoefficients}
S_{l,m}=<f(\theta,\varphi),Y_{l,m}(\theta,\varphi)>=\int_0^{2\pi}\int_0^{\pi}f(\theta,\varphi) \overline{Y_{l,m}}(\theta,\varphi) \sin\theta d\theta d\varphi.
\end{equation}
The spherical harmonic coefficients have the particularity of being relevant. Indeed, the coefficients that contain the most information are those relating to low frequencies. This property is very interesting for reconstruction since it is possible to have a reconstructed surface very close to the initial surface with a limited number of coefficients.

The most common problem in modeling with spherical harmonics is how to define the optimal reconstruction order that represents the closest model to the initial surface to be modeled. This task still remains subjective and determined always iteratively. At each iteration, the reconstruction order is incremented. The recostruction is finished when the desired model is obtained.
\subsection{New recursive spherical harmonics}
In this section we will provide a new method based on recurrence rule to show that applying spherical harmonics method for image processing did not necessitate nor the explicit or exact computation of all the spherical harmonics coefficients of the image neither the approximation yielded in the previous section. In contrast, such calculus may be issued from one level computation and sometimes from one or two coefficients. The idea is based on the theory of Legendre polynomials. These may also be introduced via the induction rule
$$
P_{n+1}(x)=\displaystyle\frac{2n+1}{n+1}xP_{n}(x)-\frac{n}{n+1}P_{n-1}(x)
$$
satisfied for all $x\in[-1,1]$, and all $n\in\mathbb{N}^{*}$ with initial data $P_{0}(x)=1$ and $P_{1}(x)=x$. Denote:
$$
a_l=\frac{2l+1}{l+1}~,~~b_l=\frac{l}{l+1}
$$
and the weight function $w_m(x)=(1-x^2)^{\frac{m}{2}}$. We know that
$$
P_{l,m}(x)=(-1)^mw_m(x)\frac{d^mP_l(x)}{dx^m}.
$$
Using the induction rule this yields that
$$
\begin{array}{lll}
\medskip\,P_{l,m}(x)&=&(-1)^ma_{l-1}w_m(x)\frac{d^m}{dx^m}(xP_{l-1})(x)\\
\medskip&&-(-1)^mb_{l-1}w_m(x)\frac{d^m}{dx^m}P_{l-2}(x).
\end{array}
$$
Denote by the next
$$
I_m=\frac{d^m}{dx^m}(xP_{l-1}(x)).
$$
Obviously,
$$
I_0=xP_{l-1}(x).
$$
For $m\geq1$, using Leibnitz rule, we get
$$
I_m=x\frac{d^m}{dx^m}P_{l-1}(x)+m\frac{d^{m-1}}{dx^{m-1}}P_{l-1}(x).
$$
Therefore,
$$
\begin{array}{lll}
\medskip\,P_{l,m}(x)&=&(-1)^mw_m(x)[a_{l-1}x\frac{d^m}{dx^m}(P_{l-1}(x))\\
\medskip&&+a_{l-1}m\frac{d^{m-1}}{dx^{m-1}}P_{l-1}(x)-b_{l-1}\frac{d^m}{dx^m}P_{l-2}(x)].
\end{array}
$$
Otherwise,
$$
\begin{array}{lll}
\medskip\,P_{l,m}(x)&=&a_{l-1}xP_{l-1,m}(x)\\
\medskip&&+ma_{l-1}P_{l-1,m-1}(x)\\
\medskip&&-b_{l-1}P_{l-2,m}(x).
\end{array}
$$
Denote by the next
$$
\alpha_{l,m}=\sqrt{\frac{(2l+1)(l-m)}{(2l-1)(l+m)}},
$$
$$
\beta_{l,m}=\sqrt{\frac{(2l+1)}{(2l-1)(l+m)(l+m-1)}}
$$
and
$$
\gamma_{l,m}=\sqrt{\frac{(2l+1)(l-1)(l-m-1)}{(2l-3)(l+m)(l+m-1)}}.
$$
We get a recurrence rule
\begin{equation}\label{SHBasisrecursive}
\begin{array}{lll}
\medskip\,Y_{l,m}&=&a_{l-1}\alpha_{l,m}\cos\theta Y_{l-1,m}\\
\medskip&&+e^{i\varphi}ma_{l-1}\beta_{l,m}\sin\theta Y_{l-1,m-1}\\
\medskip&&-b_{l-1}\gamma_{l,m}Y_{l-2,m}.
\end{array}
\end{equation}

\begin{definition}
	We call the sequence of vectors $\mathcal{F}=(a_l\;\;b_l)^T_l$ the spherical harmonics bi-filter.
\end{definition}

\section{Spherical harmonics entropy}
As we have noticed previously, entropy concept is originally related to thermodynamics. Next, Shannon's entropy has been extended or reformulated in the framework of information theory as a mathematical function that corresponds to the amount of information contained or delivered by a source of information. For a process characterized by a number N of states or classes of events, Shannon's entropy is defined by (\cite{Shannon})
\begin{equation}\label{equationdegeth0}
S=H_1=-\sum_{i=1}^{N}p_i \log(p_i)
\end{equation}
where $p_i$ is the probability of occurrence of the $i$-th event. In the context of continuous frame, Shannon's entropy is reformulated in an integral form
$$
S=-\int_0^{\infty}p(x)\log(p(x))dx,
$$
where $x$ is a point in the space domain and $p$ is the spatial probability density that guaranties the occurence of an vent at time or position $x$. 

In the literature, entropy may be expressed otherwise accordingly to the mathematical/physical basis of the space where the studied system lives such as wavelet bases. (See \cite{LabatDR}, \cite{RossoO} and \cite{LyubushinA}). An extension to 2D case has been provided in \cite{OriettaN}.
In 3D images, especially spherical ones, in our knowledge, the concept of entropy has not been yet widely known. This motivates our development in the present paper, where we intend to introduce a special expression of entropy by means of spherical harmonics' bases, which we call spherical harmonics' entropy (SHE). 
Let $f\in L^2(S^2)$ be a finite-energy image on the sphere. Observing (\ref{decompositionshseries}) and (\ref{shcoefficients}), the $l$-level $m$-detail energy of the image $f$ is defined by 
\begin{equation}\label{equationdegeth1}
E_h(l,m)=|S_{l,m}|^2.
\end{equation}
Next, as for the context of wavelet entropy, the $l$-level energy of $f$ is the sum on all $m$ details, $|m|\leq l$, and thus defined by
$$
E_{h}(l)=\frac{1}{N}\sum_{|m|\leq l}E_{h}(l,m).
$$
The total energy is the sum over all the $l$-multiresolution levels, 
$$
E_{h}=\sum_{l}E_{h}(l)=\frac{1}{N}\sum_{l}\sum_{|m|\leq l}E_{h}(l,m).
$$
In practice, of course, we could not compute all levels $l$, as these may be infinite, and thus we estimate the energy by means a finite approximation
$$
E_{h,J}=\sum_{l=0}^JE_{h}(l)=\frac{1}{N}\sum_{l=0^J}\sum_{|m|\leq l}E_{h}(l,m)
$$
which w call by the next the $J$-approximation of spherical harmonics' energy. 

Next, as for Shannon's entropy, we introduce the spherical harmonics' distribution or spherical harmonics probability density at the level $l$ by
$$
P_h(l)=\frac{E_h(l)}{E_h}=\displaystyle\frac{\displaystyle\sum_{|m|\leq l}E_{h}(l,m)}{\displaystyle\sum_{l}\displaystyle\sum_{|m|\leq l}E_{h}(l,m)}.
$$
Or equivalently, 
\begin{equation}\label{equationdegeth4}
P_{h}(l)=\displaystyle\frac{\displaystyle\sum_{|m|\leq l}|S_{i,m}|^2}{\displaystyle\sum_{l}\sum_{|m|\leq l}|S_{j,m}|^2}.
\end{equation}
Finally, the spherical harmonics' entropy is defined by 
\begin{equation}\label{equationdegeth5}
SHE=-\displaystyle\sum_{l}P_{h}(l)\log(P_{h}(l)).
\end{equation}
In practice, as previously, we compute the $J$-estimation of the SHE
\begin{equation}\label{Jequationdegeth5}
SHE(J)=-\displaystyle\sum_{l=0}^JP_{h}(l)\log(P_{h}(l)).
\end{equation}
\section{Experimentation}
In this section, some illustrative examples are developed in order to show the efficiency and the effectiveness of the theoretical concepts discussed previously. Firstly, A spherical harmonics processing is conducted in some 3D shapes to show the efficiency of such technique in 3D images processing. Next, the spherical harmonics entropy is computed to evaluate the most appropriate reconstruction order for each 3D imagee considered. The diagram bellow illustrates the steps performed for the reconstruction of the spherical harmonic model.
\newpage
\begin{figure}[h!]
	\begin{center}
		\begin{tikzpicture}
		\begin{scope}[xscale=2,yscale=1]
		\node (A1) at (0,15) [rectangle,draw] {Initial surface};
		\node (A2) at (0,13.5) [rectangle,draw] {Assesment of $(Y_0^m,Y_1^m)$ and $l=2$};
		\node (A3) at (0,12) [rectangle,draw] {Recursive construction of $Y_k^m$, $k\geq3$};
		\node (A4) at (0,10.5) [rectangle,draw] {Computation of SH coefficients $S_l^m$};
		\node (A5) at (0,9) [rectangle,draw] {Computation of SHE at order $l$};
		\node (A6) at (0,7.5) [diam] { If $SHE_l\not=0$ };
		\node (A7) at (0,6) [rectangle,draw] {Final and optimal order reached is $l+2$};
		\node (A8) at (0,4.5) [rectangle,draw] {Reconstruction of the SH model};
		\node (B1) at (2.5,12) [rectangle,draw] {$l=l+1$};
		\draw[->] (A1) -- (A2);
		\draw[->] (A2) -- (A3);
		\draw[->] (A3) -- (A4);
		\draw[->] (A4) -- (A5);
		\draw[->] (A5) -- (A6);
		\draw[->] (A6) -- node[anchor=east] {yes} (A7);
		\draw[->] (A7) -- (A8);
		\draw[<-] (A3) -- (B1);
		\draw[-] (B1) -- (2.5,7.5);
		\draw[<-] (A6) -- node[anchor=south] {no} (2.5,7.5);
		\end{scope}
		\end{tikzpicture}
	\end{center}
	\caption{Steps of reconstruction process using SHE.} \label{organigramme}
\end{figure}
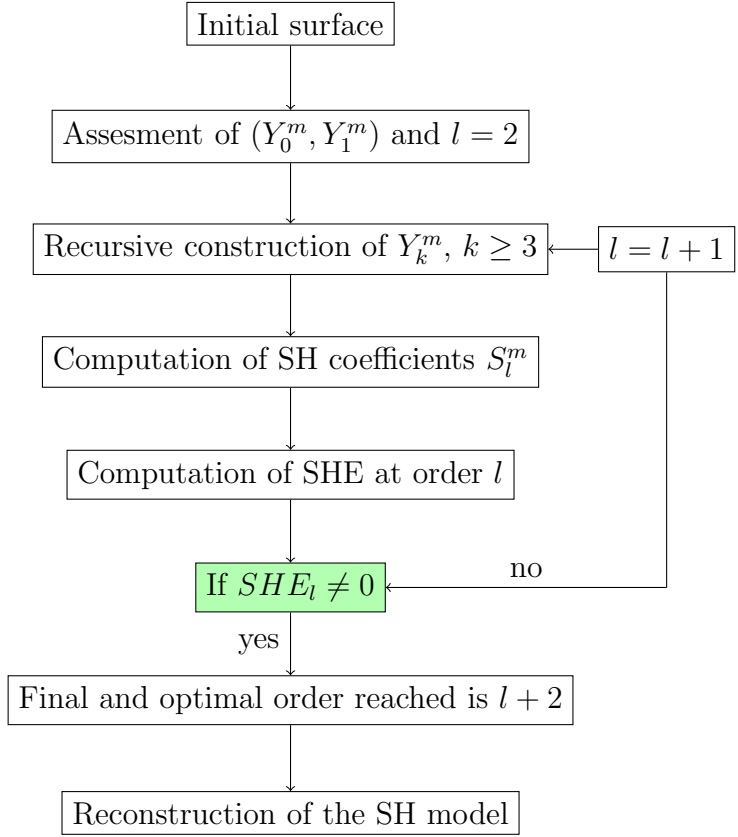
\subsection{Spherical harmonics modeling}
Spherical harmonic modeling is an extension of the powerful technique of 2D Fourier description to 3D objects. Spherical harmonics are smooth functions that can be easily arranged in increasing order of spatial frequency. Low degrees spherical harmonics capture the global shape characteristics of the surface.

$$
\begin{pmatrix}
X(\theta,\varphi)\\
Y(\theta,\varphi)\\
Z(\theta,\varphi)\end{pmatrix}
=\begin{pmatrix}\\
\displaystyle\sum_{n\in\mathbb{N}}\sum_{m=0}^{n}(a_n^m\cos(m\varphi)+b_n^m \sin(m\varphi))P_n^m(\cos(\theta))\\
\\
\end{pmatrix}\begin{pmatrix}
\cos(\theta)\cos(\varphi)\\
\cos(\theta)\sin(\varphi)\\
\sin(\theta)\end{pmatrix}.
$$

In the following figure (Fig. \ref{AnalyseFigure}), synthetic surface reconstructions are provided with different degrees or spherical harmonics multiresolution levels for a 3D shape (Fig. \ref{AnalyseFigure}. A.).
\begin{figure}[h]
	\begin{center}
		\includegraphics[scale=0.60]{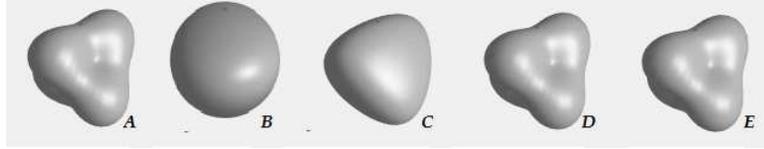}
		\caption{Reconstruction of a spherical model: A. Initial synthetic object, B. Degree 1, C. Degree 2, D. Degree 4, E. Degree 5.}
		\label{AnalyseFigure}
	\end{center}
\end{figure}\\ 
In the following illustrations, we computed for different 3D shapes the number of spherical harmonics coefficients contributing in the 'best' reconstruction relatively to the 'best' multiresolution levels estimated above.

\begin{figure}[h]
	\begin{center}
		\includegraphics[scale=0.60]{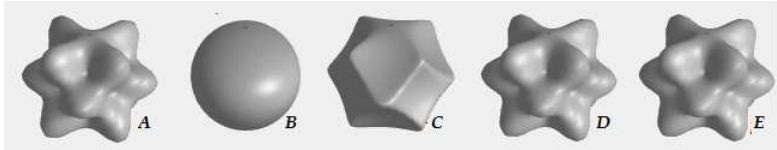}
		\caption{Reconstruction of spherical model 2: A. Initial synthetic object, B. Degree 4, C. Degree 6, D. Degree 7, E. Degree 8.}
		\label{AnalyseFigure2}
	\end{center}
\end{figure}
The second figure (Fig. \ref{AnalyseFigure2}) illustrates a second proof of efficincy of spherical harmonics processing for 3D images. It illustrates precisely a second synthetic surface reconstructions with spherical harmonics multiresolution levels for a 3D shape (Fig. \ref{AnalyseFigure2}. A.). 

We notice from these examples that for the first figure (Fig. \ref{AnalyseFigure}) we managed to rebuild the model at $l=4$. For the second model (Fig. \ref{AnalyseFigure2}) a best reconstruction was done at a multiresolution level $l=7$.\\
\begin{figure}[h]
	\begin{center}
		\includegraphics[scale=0.60]{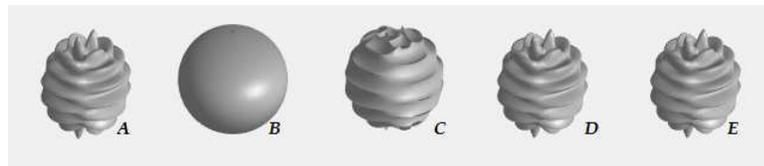}
		\caption{Reconstruction of spherical model 3:  A. Initial synthetic object, B. Degree 4, C. Degree 15, D. Degree 17, E. Degree 18.}
		\label{AnalyseFigure3}
	\end{center}
\end{figure}\\
The following and last illustrative example (Fig. \ref{AnalyseFigure3}) shows the number of contributoers coefficients in the synthetic image reconstruction at the level $l=17$.\\
\begin{table}[h]
	\begin{center}
		\begin{tabular}{||c|c|c|c||}
			\hline\hline
     		Initial surface&$N_S$ &J&$N_{HS}$\\
     		\hline
     		Fig. \ref{AnalyseFigure}&2601&4&25\\
     		\hline
     		Fig. \ref{AnalyseFigure2}&3721&7&64\\
			\hline
			Fig. \ref{AnalyseFigure3}&8281&17&324\\
			\hline\hline
		\end{tabular}
        \caption{Computational parameters}
        \label{comparisontable3}
	\end{center}
\end{table}\\
Tables above highlight the compactness characteristic of spherical harmonics and shows the fidelity of the reconstructed model compared to the initial surface. \\
\begin{itemize}
	\item $N_S$ to be the number of points of the initial surface.
	\item $J$ to be the reconstruction order or the multiresolution level.
	\item $N_{HS}$ to be the number of spherical harmonics coefficients contributors.
\end{itemize}
It is important to note that that the computation of spherical harmonic basis was performed using the new recursive spherical harmonics presented in paragraphe(3.2).This new method allows to gain in complexity. 

In fact, with the classical way, spherical harmonic basis is calculated on every point of the initial surface for every order $L$ (\ref{SHbasis}). So far an initial surface with $N$ points we will obtain a complexity of $O(N^2)$ which is important especially when complex surface is treated. 

However, with the new recursive method, spherical harmonic basis will calculated only for the $1^{st}$ order and then the following terms will be obtained with a simple deduction using (\ref{SHBasisrecursive}) we have then a $O(N)$. Which implies a considerable gain in programing and execution time.

To more illustrate the charcateristic of compactness of this modeling method, we present below the spherical harmonic coefficients spectrum for Figures \ref{AnalyseFigure}, \ref{AnalyseFigure2} and \ref{AnalyseFigure3} respectively.

The spectra above show that the main information of the 3D object lies in the low frequencies and thus, the majority of the data is centered at zero. The rest of the spectrum represents the detail of the initial surface.
\begin{figure}[h]
	\begin{center}
		\includegraphics[scale=0.50]{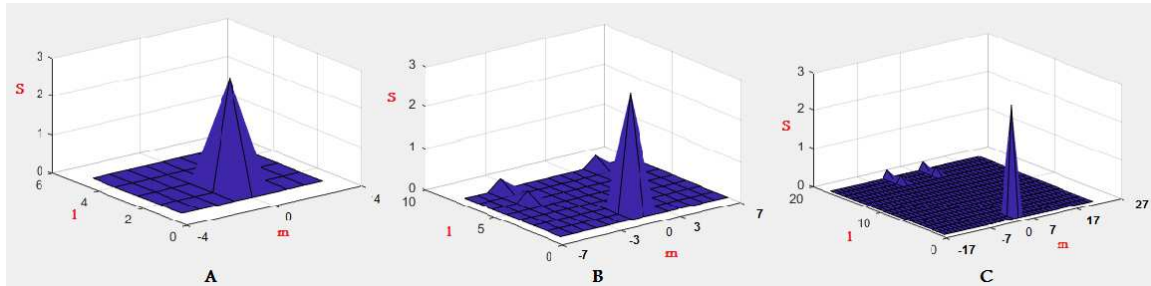}
		\caption{Spherical harmonics coefficients spectrum: A. Degree 4. B. Degree 7. C. Degree 17.}\label{SHSpectrum}
	\end{center}
\end{figure}\\

\subsection{Spherical harmonics entropy}
Remark that in all these examples illustrated above, as well as in the literature on spherical harmonics applications (in our knowledge), there is no theoretical or precisely no a priori rule that permits to predict the optimal degree or level of reconstruction. Theoretically speaking, as we know that in the Hilbert space $L^(S^2)$

$$
\|f-\sum_{k=0}^{l}\sum_{m=-k}^kS_{k,m}Y_{k,m}\|\longrightarrow0
$$
whenever $l\rightarrow+\infty$, we get more and more best levels as we increase $l$. Which makes it questionable to stop the procedure in some level.

Therefore, a natural problematic raised is how to define automatically and accurately the optimal reconstruction order. In this section, we show that spherical harmonics entropy may be a good blackbox permitting to stop the procedure with a best apprximation guaranteed.

Assessment of the spherical harmonics entropy will allow to determine in a precise way the optimal order of reconstruction. Indeed, as entropy in its general form and definition as well as mathematical/physical meaning is a type of dimension. So, it should be a somehow global measure of invariance for the studied system. Its value should therefore be stationnary or quitely constant as the multiresolution level increases. This is shown in the following part.

\begin{figure}[h]
	\begin{center}
		\includegraphics[scale=0.60]{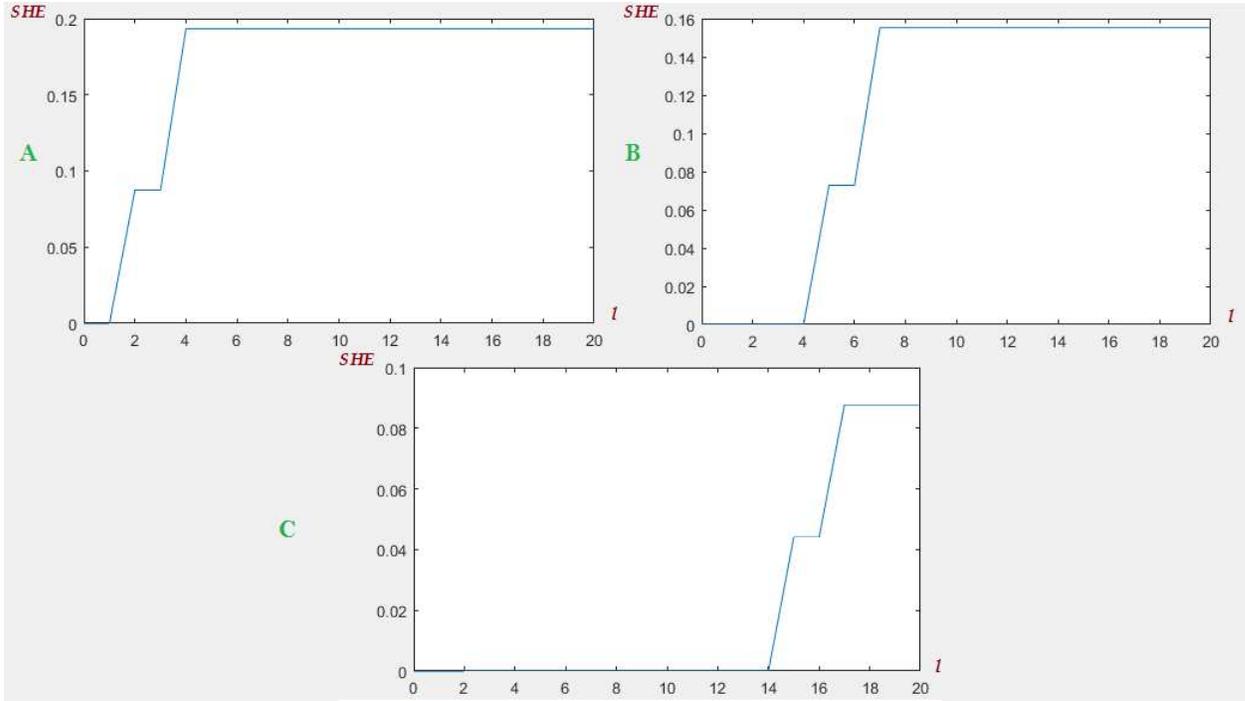}
		\caption{SHE for Fig. \ref{AnalyseFigure}, \ref{AnalyseFigure2} and \ref{AnalyseFigure3}.}
		\label{SHE-Fig4}
	\end{center}
\end{figure}

The following graphs (\ref{SHE-Fig4}) illustrates represents SHE curves calculated on three samples of genus zero surfaces shows in figures \ref{AnalyseFigure}, \ref{AnalyseFigure2} and \ref{AnalyseFigure3} respectively, during the process of spherical harmonic reconstruction. It is important to say that the sHE values is stabilized from a certain order  $L$. This order will be considered as the optimal $L$ for the reconstruction process.\\
\begin{itemize}
	\item $A$ the SHE for figure \ref{AnalyseFigure}.
	\item $B$ the SHE for figure \ref{AnalyseFigure2}.
	\item $C$ the SHE for figure \ref{AnalyseFigure3}.
\end{itemize}
It follows from these illustrations that spherical harmonics entropy is an efficient way to predict the sufficient best optimal order of reconstruction. This method allows us to gain processing time since we are no longer forced to test all the orders preceding the final order of reconstruction.

The table above shows the $SHE$ evolution when we increase the value of $L$. It is important to know that the $SHE$ stabilizes when the $SH$ model do not show more details. So, then we conclude that the optimal order $L$ is reached.
\newpage
\begin{table}[h]
	\begin{center}
	\begin{tabular}{|c|c|c|}
		\hline
		 surface& the order $L$ & Spherical harmonic entropy\\
		\hline
		\multirow{2}{*}{Fig. \ref{AnalyseFigure}}   &0   &0,0000      \\
		   	\cline{2-3}
			&1            &0,0000 \\
		    \cline{2-3}
			&2            &0,0874 \\
		    \cline{2-3}
            &3            &0,0874 \\
		    \cline{2-3}
		    &$>=4$         &0,1600 
		                   
		\\ \hline
		\multirow{2}{*}{Fig. \ref{AnalyseFigure2}}   &0..4   &0,0000      \\ 	
			\cline{2-3}
			&5           &0,0728      			\\ 	
			\cline{2-3}
			&6           &0,0728      			\\ 	
			\cline{2-3}
			&$>=7$          &0,1554     
		\\ \hline
		\multirow{2}{*}{Fig. \ref{AnalyseFigure3}}   &0..14   & 0,0000    \\ 	
 			\cline{2-3}
			&15           &0,0440 \\ 	
			\cline{2-3}
			&16           &0,0440 \\ 	
			\cline{2-3}
			&$>=17$           &0,0875
		\\ \hline
	\end{tabular}
		\caption{Computational parameters for Fig. \ref{SHE-Fig4}.}\label{comparisontable4}
	\end{center}
\end{table}
Experimentation is performed for a lot of number of surfaces of varying complexity. Result avec always the same. This is why we conclude that the $SHE$ may be a precise way to know from the begining depending of spherical harmonics coefficients the optimal order of spherical harmonic model.\\

\section{Conclusion}
In this paper, spherical harmonics have been examined and proved to be efficient candidates in the reconstruction of spherical image processing. In the present work, the stimulating idea is based on the connection between well-known spherical harmonic coefficients and signal energy. Then, we introduced the notion of spherical harmonic entropy (SHE) in the reconstruction of the spherical model of a 3D object. Fast and accurate algorithms have been obtained. moreover the new definition of entropy 3d based on the coefficients of spherical harmonics has improved the work, it makes it possible to measure the complexity of model to be reconstructed, easily and quickly, and subsequently to have a perfect reconstruction in real time and without loss of information.
For our future work, we intend to validate the proposed method on more complex surfaces represented under triangular mesh. To do this we must provide a phase of spherical parametrization which consists of projecting the initial object on the unit sphere before calculating the spherical harmonic coefficients.\\

\end{document}